\documentclass[10pt]{article}
\DeclareMathAlphabet{\mathpzc}{OT1}{pzc}{m}{it}

\usepackage{latexsym}
\usepackage{marvosym}
\usepackage{amssymb}
\usepackage{amsfonts}
\usepackage{amsmath}
\usepackage{amsthm}
\usepackage{dsfont}    
\usepackage{bbold}     
\usepackage[english]{babel}
\usepackage{caption}
\usepackage{epsfig}
\usepackage{float}
\usepackage{subfigure}
\usepackage{psfrag}
\usepackage{graphicx}
\usepackage{epsfig}
\usepackage{relsize}
\usepackage{comment}
\usepackage{enumerate}

\newcommand{\ttri}{\mathpzc{T}_z^{\mathrm{TRI}}}
\newcommand{\ttide}{\mathpzc{T}_z^{\mathrm{TIDE}}}
\newcommand{\atri}{\ddot\theta^{\,\mathrm{TRI}}}
\newcommand{\atide}{\ddot\theta^{\,\mathrm{TIDE}}}
\newcommand{\Qtri}{\mathpzc{Q}^{\mathrm{TRI}}}
\newcommand{\Qtid}{\mathpzc{Q}^{\mathrm{TIDE}}}
\newcommand{\sgn}{\mathop{\mathrm{sgn}}}
\newcommand{\mtd}{\langle\dot\theta\rangle}
\newcommand{\mtdn}{\langle\dot\theta/n\rangle}

\def\v#1{\mathbf{#1}}

\newcommand{\ymax}{\dot\theta_{\mathrm{max}}}
\newcommand{\ymin}{\dot\theta_{\mathrm{min}}}

\begin{document}
\title{Numerics for the spin orbit equation of Makarov\\
with constant eccentricity}

\author
{\bf Michele Bartuccelli$^1$, Jonathan Deane$^1$, Guido Gentile$^2$
\vspace{2mm}
\\ \small 
$^1$ Department of Mathematics, University of Surrey, Guildford, GU2 7XH,
UK 
\\ \small
$^2$ Dipartimento di Matematica e Fisica, Universit\`a di Roma Tre, Roma, I-00146,
Italy
\\ \small 
E-mail: m.bartuccelli@surrey.ac.uk, j.deane@surrey.ac.uk,
gentile@mat.uniroma3.it
}
\date{}
\maketitle

\begin{abstract}
We present an algorithm for the rapid numerical integration of a
time-periodic ODE with a small dissipation term that is $C^1$ in the velocity.
Such an ODE arises as a model of spin-orbit coupling in a star/planet
system, and the motivation for devising a fast algorithm for its solution
comes from the desire to estimate probability of capture in various
solutions, via Monte Carlo simulation: the integration times are very long,
since we are interested in phenomena occurring on times similar to the
formation time of the planets.  The proposed algorithm is based on 
the High-order Euler Method (HEM)
which was described in~\cite{hem}, and it requires computer algebra to set up the
code for its implementation. The pay-off is an overall increase in speed by a
factor of about $7.5$ compared to standard numerical methods.
Means for accelerating the purely numerical computation are also discussed.
\end{abstract}

\section{Introduction}
\label{intro}

The paper first concisely sets out the model for spin-orbit coupling in, for instance,
the Sun-Mercury system, as described by Makarov~\cite{makarov} and further discussed by
Noyelles, Frouard, Makarov and Efroimsky in~\cite{nfme}.  We refer to this model as
the NFME model throughout. We then describe a fast numerical algorithm, based on the
High-order Euler Method (HEM) introduced in~\cite{hem}, for solving the spin-orbit ODE
incorporating the NFME model. After rigorous testing, we conclude that the
proposed algorithm is about 7.5 times faster than simply using a
general-purpose numerical ODE solver.

There are two main reasons for developing a fast algorithm for solving this
version of the spin-orbit ODE, bearing in mind that the ultimate objective
is to establish, by Monte Carlo simulation, the probability of capture in the various 
attracting solutions of the ODE. The first reason is that, since the
dissipation in the problem is low, very long integration
times --- of the order of that of the formation of the planets --- are
needed in order to establish capture. The second reason is that to estimate
capture probabilities with high confidence, solutions starting from many
randomly-chosen initial conditions need to be
considered: in fact, the width of the confidence interval decreases only as
the square root of the number of solutions investigated. In most of the previous
related work, for instance, \cite{CL}, \cite{CC1}, \cite{makarov},
\cite{makarov2012}, only about 1000 initial conditions were considered. By
contrast, we are able to consider more than 50,000 here, for which the total computation
time is of order two weeks.

The dynamics of the solutions of the spin-orbit ODE with the NFME model are also
interesting, and we discuss these in a separate paper~\cite{paper2}.

In~\cite{nfme}, the authors point out how the model that they advocate describes more
realistically the effect of tidal dissipation, in comparison with the much simpler constant 
time-lag (CTL) model. The latter has been used extensively in the
literature since the seminal paper by Goldreich and Peale \cite{GP} --- see
for instance~\cite{CL}, \cite{CC1}, \cite{CC2} and \cite{hem} for recent
results on the basins of attraction of the principal attractors. A
clear drawback of the CTL model is that, for the parameters of the Sun-Mercury
system, the main attractor (with 70\% probability of being observed) is quasi-periodic,
with an approximate frequency of 1.256. In fact, Mercury is in a 3:2 spin-orbit resonance,
and the CTL model indicates that the probability of this is only about 8\%.
On the other hand, with the NFME model, all attractors have mean velocity
close to rational values, and about 42\% of the initial data are captured
by an attractor with mean velocity close to 3/2.

The greater realism of the NFME model comes at a cost, however:

\begin{enumerate}
\item The NFME model is not smooth; in fact, it is only $C^1$ in the
angular velocity, $\dot\theta$, as we shall see.
\item There is considerable detail to take into account in implementing the
NFME model mathematically.
\item The functions in the model take considerably more computation time to
evaluate, and so numerical calculations are comparatively slow.
\end{enumerate}

The first point is important both for application of perturbation theory to
the problem as well as for implementing an appropriate fast HEM-type algorithm to
solve the ODE. Both of these objectives are important. For
instance, perturbation theory is a useful tool for obtaining any kind of analytical
results for the dynamics; and the motivation for faster algorithms has already been
discussed above.

Throughout this paper, the units for mass, length and time will be kg, km
and Earth years, yr, respectively. For ease of cross-referencing,
an equation numbered (N.$x$) in this paper is equation ($x$) in~\cite{nfme},
modified if necessary. When we refer to a `standard numerical method' or `numerical
method' for ODE solving, we mean one of the family of well-known general purpose ODE
solvers (e.g. Runge-Kutta, Adams); when we mean HEM, which is of course
also numerical, we refer to it explicitly by name.

\section{The NFME model}

\subsection{The spin-orbit ODE}

The ODE is
\begin{equation}
\tag{N.1}
\ddot\theta = \frac{\ttri(\theta, t) + \ttide(\dot{\theta})} {\xi M_{planet} R^2}
\label{nfme_ode}
\end{equation}
where $\xi$ is a measure of the inhomogeneity of the planet (with $\xi = 2/5$
for a homogeneous sphere), $M_{planet}$ is the mass of the planet and $R$ is its radius. 
Numerical values for these, and indeed all relevant parameters can be found in
Table~\ref{numvals}.

\subsection{The triaxiality torque}

The triaxiality-caused torque along the $z$-axis, $\ttri(\theta, t)$, is a torque 
exerted on the planet by the gravitational field of the star, arising from
the fact that the planet is not a perfect sphere.
Goldreich and Peale~\cite{GP} and NFME approximate $\ttri$ by its quadrupole part,
which is given by
\begin{equation}
\tag{N.5}
\ttri \approx -\frac{3a^3}{2r^3}(B-A) n^2\sin 2(\theta - f)
= -\frac{3a^3}{2r^3}(B-A) n^2 (\sin 2\theta\cos 2f - \cos 2\theta\sin 2f)
\label{tri_def}
\end{equation}
where the principal moments of inertia of the planet, along the $x$, $y$
and $z$ axes respectively, are $A < B < C = \xi M_{planet} R^2$; $n$
is its mean motion; $r = r(t)$ is the distance between the centres of mass of the
star and the planet at time $t$; $a$ is the semi-major axis of the orbit of the planet;
$\theta(t)$ is the sidereal angle of the planet, measured from the major axis of its orbit;
and $f = f(t)$ is the true anomaly as seen from the star. We further define $\mathpzc{M}(t)$,
the mean anomaly, which is such that $n = \dot{\mathpzc{M}}$, with a good
approximation~\cite{nfme} being given by\footnote{Numerically, this is true:
\cite{nfme} give $n = 26.0879$ yr$^{-1}$ for Mercury and the approximation yields
$n\approx 26.0897$ yr$^{-1}$.}
$n \approx \sqrt{G(M_{planet} + M_{star})/a^3}$, where $G$ is the
gravitational constant and $M_{star}$ is the mass of the
star. Hence, aside from a constant of integration which we set to zero, $\mathpzc{M} = nt$.

The standard procedure from this point is to make the following pair of Fourier
expansions:
\begin{equation}
\tag{N.7, N.8}
\left(\frac{r(t)}{a}\right)^n \begin{array}{c}\cos\\ \sin\end{array} \left(m f(t)\right)
= \sum_{k\in\mathbb Z} X_k^{n, m}(e) \begin{array}{c}\cos\\ \sin\end{array} \left(k n t\right)
\end{equation}
where $X_k^{n, m}(e)$ are the Hansen coefficients, which depend on the
orbital eccentricity $e$. NFME compute them via~\cite{duriez}
\begin{equation}
\tag{N.9}
X_k^{n, m}(e) = (1 + z^2)^{-n-1}\sum_{g = 0}^\infty 
(-z)^g\sum_{h = 0}^g C^{n, m}_{g-h, h}\, J_{k-m+g-2h}(ke),
\label{hancomp}
\end{equation}
where $z = (1 - \sqrt{1-e^2})/e$, $J_k(x)$ is the $k$-th order Bessel
function of the first kind, and
$$C^{n, m}_{r, s} = \left(\begin{array}{c}n+1+m\\r\end{array}\right)
\left(\begin{array}{c}n+1-m\\s\end{array}\right)$$
with the binomial coefficients (extended to all integer arguments) being given by
$$\left(\begin{array}{c}l\\m\end{array}\right) =
\begin{cases}
\frac{l(l-1)\ldots (l-m+1)}{m!} & l\in\mathbb Z, m\geq 1\\
1 & l\in\mathbb Z, m = 0\\
0 & l\in\mathbb Z, m < 0.\end{cases}$$

For convenience, we define the functions $G_{lpq}(e) = X_{l-p+q}^{-l-1,\,l-2p}$.
Then finally,
\begin{equation}
\tag{N.10}
\ttri(\theta, t) = -\frac{3}{2} (B-A) n^2\sum_{q\in\mathbb Z} G_{20q}(e)
\sin (2\theta - (q+2)nt).
\label{ttri}
\end{equation}
In practice, NFME sum over $q\in\Qtri = \{-4, \ldots, 6\}$.

\subsection{The tidal torque}

NFME forcefully point out the problems associated with the constant
time-lag (CTL) model, which leads to a dissipation term of the form $\gamma(\dot\theta
- \omega)$, with $\gamma$ and $\omega$ constants. This model, among other things, can give
rise to quasi-periodic solutions for the spin-orbit problem~\cite{CL,CC2,hem},
and the fact that these have not been observed in practice is evidence that the CTL model is in
some way unphysical.

NFME model tidal dissipation by expanding both the tide-raising potential
of the star, and the tidal potential of the planet, as Fourier series. The
Fourier modes are written $\omega_{lmpq}$, where $l, m, p, q$ are integers.
Appendix B of~\cite{nfme} justifies the simplification $l = m = 2$, $p = 0$,
which arises from the smallness of the obliquity of the orbit
(axial tilt), in combination with an averaging argument. This then gives the
following approximation for the polar component of the tidal torque:
\begin{equation}
\tag{N.11a}
\ttide(\dot\theta) = \frac{3 G M_{star}^2}{2 a} \left(\frac{R}{a}\right)^5
\sum_{q, j\in\Qtid} G_{20q}(e) G_{20j}(e) Q_2(\omega_{220q})
\sgn(\omega_{220q}) \cos[(q-j)\mathpzc{M}].
\end{equation}
Here, $\Qtid = \{-1,\ldots, 7\}$; 
$Q_2(\omega_{220q}) = k_2(\omega_{220q})\sin\left|\epsilon_2(\omega_{220q})\right|$,
where both $k_2$ (a positive definite, even function of the mode $\omega_{lmpq}$) and
$\epsilon_2$ (an odd function of the mode) are functions to be given shortly;
and $\sgn()$ is the usual signum function. Hence, overall, $Q_2(x)\sgn(x)$
is an \textit{odd} function of $x$.  The reason why $Q_2$ is written in this form is that
it will turn out to depend on a \textit{fractional power} of $|\omega_{lmpq}|$: the
oddness of $Q_2(x)\sgn(x)$ is thereby retained without raising a negative number to this
fractional power.

The additive error terms for this expression are $O(e^8
\epsilon)$, $O(i^2\epsilon)$ and $O\left(\epsilon R^7/a^7\right)$, in which
$i$ is the obliquity of the orbit and $\epsilon$, which is $O(1)$, is a phase lag.
For Mercury, $e^8\approx 10^{-8}$, $i^2\approx 10^{-7}$ and
$(R/a)^7 \approx 10^{-31}$.  NFME then make the additional approximation that,
since the terms with $q\neq j$ oscillate, they naturally affect the
detailed dynamics of the planet, but nonetheless, the
overall capture probabilities are insensitive to them~\cite{makarov2012}.
Hence, the tidal torque is well approximated by the secular part only,
i.e.\ the $q = j$ terms; this gives
\begin{equation}
\tag{N.11b}
\ttide(\dot\theta) = \langle\ttide(\dot\theta)\rangle_{l=2} =
\frac{3 G M_{star}^2}{2 a} \left(\frac{R}{a}\right)^5
\sum_{q\in\Qtid} G^2_{20q}(e)\, Q_2(\omega_{220q}) \sgn(\omega_{220q}).
\label{sec_ttid}
\end{equation}
Neglecting precessions, we have 
\begin{equation}
\tag{N.12}
\omega_{lmpq} = (l - 2p + q)n - m\dot\theta
\end{equation}
with $l\geq 2$, $m, p = 0, \ldots, l$ and $q\in\mathbb Z$. For the reasons
mentioned earlier, we consider only the mode
$\omega_{220q} = (q+2)n - 2\dot\theta$.

Finally, we give the NFME model for the quality function, $Q_2(\omega_{lmpq})$.
We use the abbreviation $\chi = \left|\omega_{lmpq}\right|$, in terms of
which 
\begin{equation}
\tag{N.15 mod.}
Q_2(\omega_{lmpq}) = -\frac{3 A_l}{2(l-1)}\frac{\mathpzc{I}'(\chi)\chi}
{(\mathpzc{R}'(\chi) + A_l\chi)^2 + \mathpzc{I}'(\chi)^2},
\label{Q2}
\end{equation}
where 
$$A_l = \frac{4\pi(2l^2 + 4l + 3)\mu R^4}{3l G M_{planet}^2},$$
with $\mu$ being the unrelaxed rigidity; 
\begin{equation}
\tag{N.16 mod.}
\mathpzc{R}'(\chi) = \chi + \chi^{1-\alpha}\tau_A^{-\alpha}\cos(\alpha\pi/2)\,\Gamma(\alpha + 1)
\end{equation}
and
\begin{equation}
\tag{N.17 mod.}
\mathpzc{I}'(\chi) = -\tau_M^{-1} - 
\chi^{1-\alpha}\tau_A^{-\alpha}\sin(\alpha\pi/2)\,\Gamma(\alpha + 1).
\label{RpIp}
\end{equation}
In these expressions, $\Gamma$ is the usual gamma function,
$\tau_A$ and $\tau_M$ are the Andrade and Maxwell
times of the mantle respectively and $\alpha$ is the Andrade parameter.

\subsection{Parameter values for the Sun/Mercury system}
\label{parvals}

In Table~\ref{numvals} we give numerical values appropriate to the Sun and Mercury
for the parameters introduced so far. We standardise the units to kg for
mass, km for length and (Earth) years for time.

\begin{table}[h]
\centering
\begin{tabular}{lll} \hline
\multicolumn{3}{c} {\raisebox{2.2ex}{ }Parameter values specific to Mercury}\\ \hline
Name & Symbol & Numerical value\\ \hline
Semimajor axis 		& $a$ 	& $5.791\times 10^7$ km\\
Mean motion		& $n$ 	& 26.0879 rad yr$^{-1}$\\
Planetary radius 	& $R$	& $2.44\times 10^3$ km\\
Dimensionless m.o.i. (inhomogeneity) & $\xi = C/(M_{planet} R^2)$ & 0.346\\
Triaxiality 		& $(B - A)/C$ & $9.350\times 10^{-5}$\\
Planetary mass 		& $M_{planet}$ & $3.301\times 10^{23}$ kg\\
Unrelaxed rigidity	& $\mu$ & $7.967\times 10^{28}$ kg km$^{-1}$ yr$^{-2}$\\
Present-day orbital eccentricity & $e$ & $0.2056$\\
Andrade time	& $\tau_A$ & $500$ yr\\
Maxwell time	& $\tau_M$ & $500$ yr\\
Andrade parameter	& $\alpha$ & $0.2$\\ \hline
\multicolumn{3}{c} {\raisebox{2.2ex}{ }Acceleration constants}\\ \hline
Triaxiality acceleration constant	& $\zeta$ & $0.09545$ yr$^{-2}$\\
Tidal acceleration constant	& $\eta$ & $0.03096$ yr$^{-2}$\\ \hline
\multicolumn{3}{c} {\raisebox{2.2ex}{ }Other parameter values}\\ \hline
Mass of Sun 		& $M_{star}$ & $1.989\times 10^{30}$ kg\\
Gravitational constant 	& $G$ & $6.646\times 10^{-5}$ kg$^{-1}$ km$^3$ yr$^{-2}$\\
Triaxiality index range	& $\Qtri$ & $\{-4, \ldots, 6\}$\\
Tidal index range	& $\Qtid$ & $\{-1, \ldots, 7\}$\\ \hline
\end{tabular}
\caption{Numerical values for the parameters used in the NFME model for the Sun and Mercury.}
\label{numvals}
\end{table}

\subsection{The triaxiality and tidal angular accelerations}

In order to study the spin-orbit problem further, it is convenient to
introduce two quantities with the dimensions of angular acceleration
$$\atri(\theta, t) = \ttri(\theta, t)/\xi M_{planet}R^2
\;\;\;\mbox{ and }\;\;\;
\atide(\dot\theta) = \ttide(\dot\theta)/\xi M_{planet}R^2.$$
Taken individually, these quantities are not accelerations \textit{per se}
--- more correctly, they might be described as the triaxiality/tidal
contributions to the angular acceleration, since, when computed on the
solution $\theta(t)$, their sum is
$\ddot\theta$, as we shall see from equation~(\ref{ode}) --- but for
brevity, we refer to them as the triaxiality and tidal accelerations. We then define 
\begin{equation}
\label{ddttridef}
\zeta = \frac{3(B-A)n^2}{2\xi M_{planet} R^2} \approx 0.09545\mbox{ yr}^{-2}
\;\;\;\mbox{ so that }\;\;\;
\atri(\theta, t) = 
-\zeta\sum_{q\in\Qtri} G_{20q}(e) \sin (2\theta - (q+2)nt),
\end{equation}
and
$$\eta = \frac{3G M_{star}^2 R^5}{2 a^6\xi M_{planet} R^2}\cdot\frac{3A_l}{2(l-1)}
= \frac{3\pi(2l^2 + 4l+3)}{l(l-1)}\frac{\mu M_{star}^2 R^7}{\xi M_{planet}^3 a^6}
\approx 0.03096\mbox{ yr}^{-2}$$
for $l = 2$, so that
\begin{equation}
\atide(\dot\theta) = 
-\eta\sum_{q\in\Qtid} G^2_{20q}(e) P_2\left(\left|\omega_{220q}\right|\right)
\sgn(\omega_{220q})
\label{ddttid}
\end{equation}
where 
$$P_2(\chi) = \frac{\mathpzc{I}'(\chi)\chi}
{(\mathpzc{R}'(\chi) + A_2\chi)^2 + \mathpzc{I}'(\chi)^2}$$
and $\mathpzc{R}'(\chi)$ and $\mathpzc{I}'(\chi)$ are defined in equation~(\ref{RpIp}).

In terms of these angular accelerations, the spin-orbit ODE becomes
\begin{equation}
\ddot\theta + \zeta\sum_{q\in\Qtri} G_{20q}(e) \sin (2\theta - (q+2)nt)
+ \eta\sum_{q\in\Qtid} G^2_{20q}(e) P_2\left(\left|2\dot\theta - (q+2)n\right|\right) = 0
\label{ode}
\end{equation}
and this is our starting point for the rest of the paper.

\section{The NFME model with practical values}

We now look at the behaviour and magnitude of the various functions that go
to make up the NFME model, using the values for the Sun and Mercury given in
the Table~\ref{numvals}.

\subsection{Hansen coefficients}

\begin{figure}[!ht]
\centering 
\includegraphics*[width=3in]{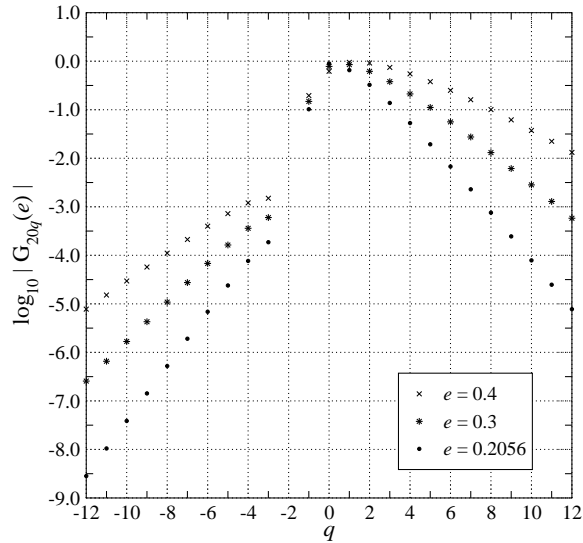}
\caption{Plot of $\log_{10} \left| G_{20q}(e)\right|$ for $e = 0.2056,\, 0.3,\, 0.4$
(the latter two values being for illustration) and
$q = -12, \ldots, 12$. Note that, for all values of $e$ considered, (i) only 
$G_{2,0,-1}$ is negative, and (ii) $G_{2,0,-2} = 0$ and is therefore not plotted. }
\label{hansen}
\end{figure} 

In figure~\ref{hansen} we plot the Hansen coefficients relevant to the
problem, $G_{20q} = X_{q+2}^{-3, 2}$, for $e = 0.2056, 0.3$ and 0.4 and for
$q = -12,\ldots, 12$. These were computed from equation~(\ref{hancomp}) by
truncating the first sum at $p = 120$ and using 20 significant figures for
computation --- these values allow the computation of the coefficients to
a much higher precision than necessary.

\subsection{The triaxiality acceleration}

We first make a simple estimate of the order of magnitude of $\atri(\theta, t)$.
Treating $\theta$ and $t$ as independent variables, it is clear from
equation~(\ref{ttri}) that $\atide\in [-D, D]$ for
all $\theta$, $t$, where
\begin{equation}
D(e) = \zeta\sum_{q\in\Qtri} \left|G_{20q}(e)\right|.
\label{Dest}
\end{equation}
We find $D(0.2056) = 0.2096$, $D(0.3) = 0.3016$ and $D(0.4) = 0.4396$ yr$^{-2}$.

\begin{figure}[!ht]
\centering 
\includegraphics*[width=3in]{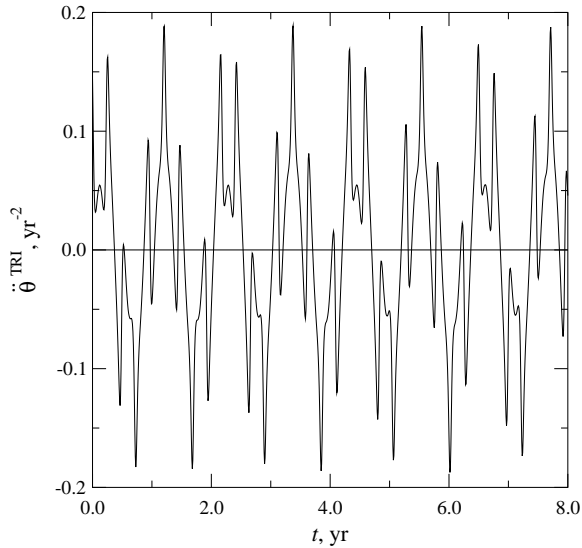}
\caption{An illustrative plot of the triaxiality angular acceleration,
$\atri(\dot\theta)$ as a function of time, with $e = 0.2056$, and $\theta(t) = 2 + 29t$.}
\label{triax}
\end{figure}

In figure~\ref{triax} we give a plot of $\atide(\theta_0 + \dot\theta_0 t, t)$
for $e = 0.2056$, $\theta_0 = 2.0$, $\dot\theta_0 = 29.0$, the latter two values being
(almost arbitrarily) chosen to approximate a solution starting from an angular
velocity slightly greater than $n$. Note that $\atide\in [-0.19, 0.19]$, which is
consistent with the estimate, from equation~(\ref{Dest}),
of $[-0.21, 0.21]$ when $e = 0.2056$.

\subsection{The tidal acceleration}

\begin{figure}[!p]
\centering 
\includegraphics*[width=3in]{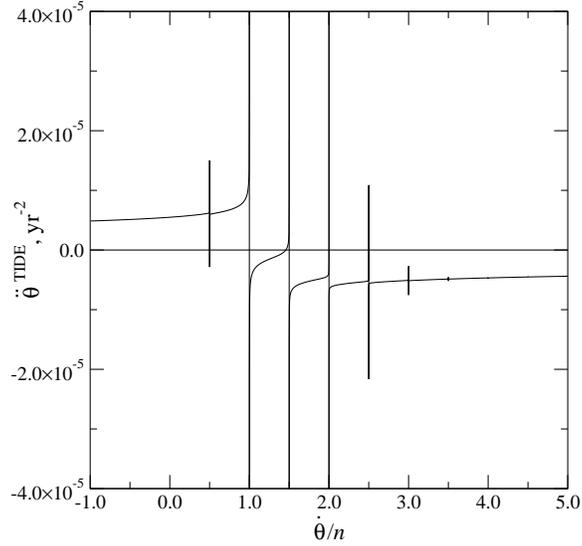}
\caption{The tidal angular acceleration, $\atide(\dot\theta)$,
with $e = 0.2056$, as given by equation~(\ref{ddttid}).
`Kinks' occur at values $\dot\theta/n = 1/2, 1, 3/2, 2, \ldots, 9/2$,
although those at $4, 9/2$ are too small to see on this scale. Note that the
angular acceleration does not change sign at a kink for $\dot\theta/n = 3, 7/2, 4, 9/2$.}
\label{kink1}
\end{figure}

\begin{figure}[!p]
\centering 
\includegraphics*[width=3in]{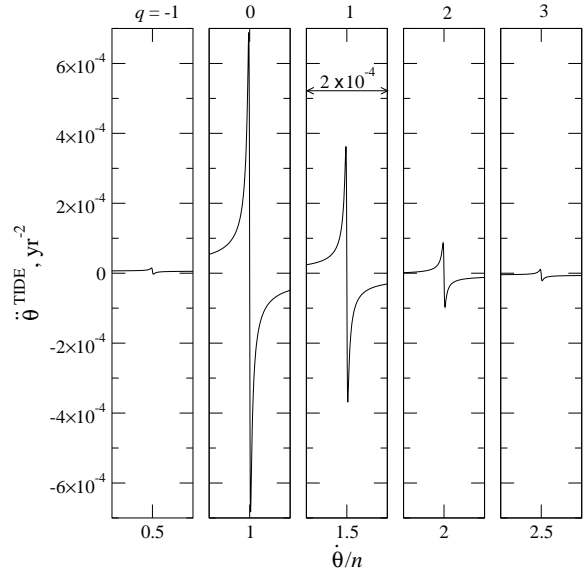}
\caption{Magnified version of the `kinks' corresponding to $q = -1, 0, 1,
2, 3$, as shown in figure~\ref{kink1}. The width of each plot is $2\times
10^{-4}$ and the vertical range is $\pm 7\times 10^{-4}$ yr$^{-2}$.}
\label{kink2}
\end{figure}

\begin{figure}[htbp]
\centering 
\includegraphics*[width=3in]{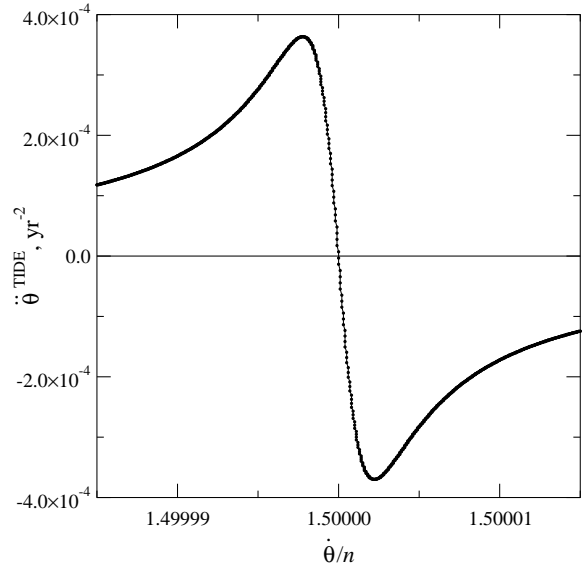}
\caption{Further magnified version of the `kink' around $\dot\theta/n = 3/2$,
suggesting that the derivative is finite for all values of $\dot\theta/n$.}
\label{kink3}
\end{figure}

We now discuss the tidal acceleration $\atide(\dot\theta)$.
Throughout this section, we take $e = 0.2056$.
Figures~\ref{kink1}--\ref{kink3} show the tidal acceleration, 
plotted vertically, versus the relative rate of rotation, $\dot\theta/n$.

This function is illustrated in Figure~\ref{kink1}, which shows
$\atide(\dot\theta)$, plotted
over the entire range of interest, $\dot\theta/n\in [-1, 5]$. The dominant 
features here are the `kinks', which occur at 
$\dot\theta/n = 1 + q/2$ for $q = -1, \ldots, 7$, this being the range of the sum defining
the tidal acceleration.

The five `kinks' at which $\atide$ changes sign are those
corresponding to $q = -1,\ldots 3$, and, for the purpose of comparison,
these are plotted over the same
narrow range of $\dot\theta/n = 1 + q/2\pm 10^{-4}$, and on the same vertical
scale, $\pm 7\times 10^{-4}$, in figure~\ref{kink2}.

\begin{figure}[!ht]
\centering 
\includegraphics*[width=3in]{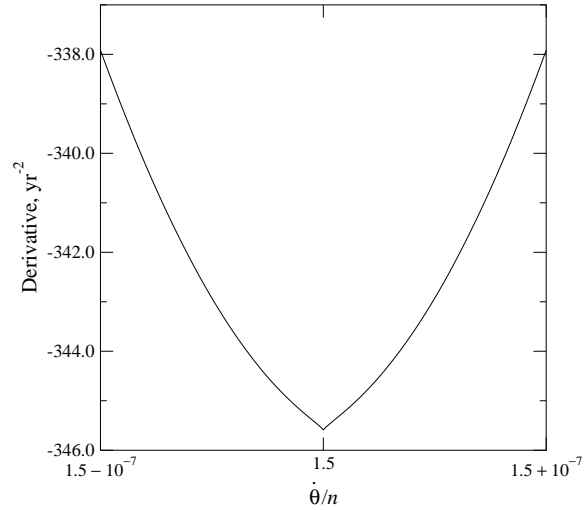}
\caption{The derivative of the tidal angular acceleration with respect to the relative
rate of rotation, that is, $n\,$d$\atide$/d$\dot\theta$, showing the
cusp at $\dot\theta/n = 3/2$.}
\label{deriv}
\end{figure}

The fact that $\atide$ is continuous, even at a
kink, is implied in figure~\ref{kink3}, in which $\atide$
has been plotted against $\dot\theta/n$ for $\dot\theta/n \in [3/2-4\times
10^{-5},\; 3/2 + 4\times 10^{-5}]$. However, the derivative
$n\,\mathrm{d}\atide/\mathrm{d}\dot\theta$ displays a
cusp at the values of $\dot\theta/n$ corresponding to a kink, as is
apparent from the definition of $Q_2(\omega_{lmpq})$ in equation~(\ref{Q2}),
in particular from its dependence on $\chi = |\omega_{220q}| = |2\dot\theta - (q+2)n|$.
This is illustrated in figure~\ref{deriv}.

\section{A fast algorithm}

We find several different solutions to the differential equation~(\ref{ode}), which one
is observed depending on the initial conditions. As well as periodic solutions, some of the
solutions we find appear, numerically at least, not to be simply periodic --- we shall have
more to say about these solutions in~\cite{paper2}. Here, we just point out that
we distinguish the different captured solutions from each other solely by their
mean $\dot\theta$ values, $\mtd$, defined by
$\mtd  = \lim_{T\rightarrow\infty}T^{-1}\int_0^T \dot\theta(t)dt$;
and that capture can take place in solutions for which a term in the sum~(\ref{ttri}), the
triaxiality torque, is
approximately zero, \textit{i.e.}\ for $\theta \approx (q+2)nt/2$, $q\in\Qtri$. Hence, we
only expect solutions for which $\mtdn \approx -1, -1/2, 1/2,\ldots, 7/2$
and $4$.

In order to compute probabilities of capture by the different solutions, with small confidence
intervals, via Monte Carlo simulation, many solutions to equation~(\ref{ode}),
starting from random initial conditions in the set 
$\mathcal Q = [0, \pi]\times [\ymin, \ymax]$, must be computed. Many are
needed because,
for a given confidence level, the width of the confidence interval is
proportional to $I^{-1/2}$, where $I$ is the number of solutions computed:
more simulations narrow the confidence interval, but rather slowly.

We take $\mathcal Q = [0, \pi]\times [0, 5n]$ in what follows. The
triaxiality acceleration depends on $2\theta$ and so we need only
consider $\theta_0\in [0, \pi]$; and the right-most kink is at $\dot\theta
= 9n/2$, so we choose the maximum value of $\dot\theta_0$ to be somewhat greater than this.

We now describe a fast algorithm for solving equation~(\ref{ode}), which is
based on the high-order Euler method (HEM), described in detail in~\cite{hem}.
This method has to be adapted for the NFME model because of the
discontinuities in the second derivative of $\atide(\dot\theta)$ --- that is, at the centres of the
kinks --- that occur at $\dot\theta/n = 1/2,
1, 3/2, \ldots, 9/2$.  The adaptation we use requires that the subset
$\mathcal Q$ of the $(\theta, \dot\theta)$ phase plane be split into strips,
all with $\theta\in [0, \pi)$, but with a variety of ranges of $\dot\theta$.
This splitting is needed in order to meet the error criterion --- see below.
A strip can be of type H (HEM), when it is far enough removed
from the kinks that the HEM can be used; and type N (Numerical), surrounding
a kink, where a suitable numerical method has to be used.
In a type N strip, there is a possibility that the ODE is stiff, and 
the numerical method chosen should take this into account.

\begin{figure}[!ht]
\centering 
\includegraphics*[width=5in]{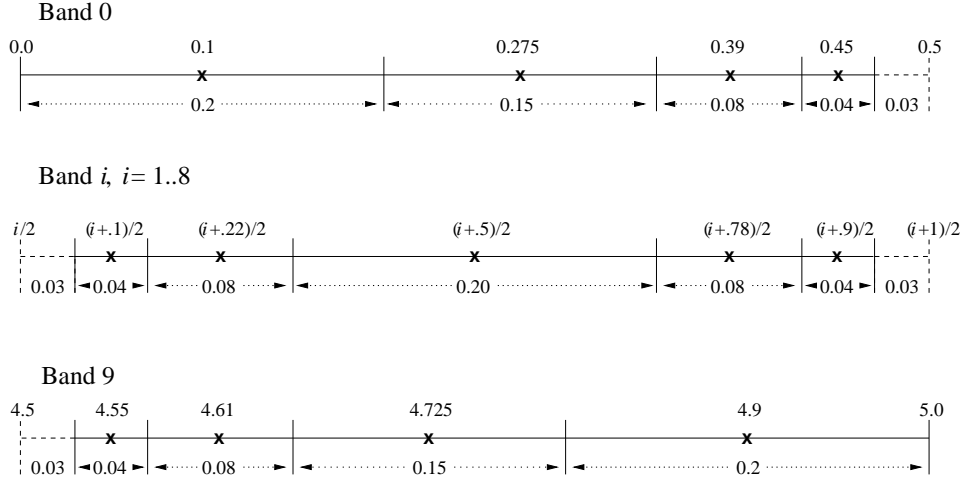}
\caption{The different strips used for H (HEM, solid line) and 
N (numerical, dashed line) computation. Numerical values of $\dot\theta/n$ are shown. A
`{\sf x}' shows the expansion point for a series solution. The width of a
region in terms of $\dot\theta/n$ over which the series solution is
valid, for all $\theta$, to within the error bounds (see text) is shown below the line.
The strips are grouped into bands. In Bands 0 and 9 there are four H strips and one N 
strip; in the remaining bands, there are five H, and two N strips.}
\label{bands}
\end{figure} 

The detail of this splitting is given in figure~\ref{bands}, in which
$\dot\theta/n$ is plotted horizontally, with type H regions being shown as
continuous lines, type N as dashed lines. The strips are grouped into ten bands,
one band covering the region between two kinks. For instance, Band 0 consists
of four type H strips, centred on $\dot\theta/n = 0.1, 0.275, 0.39$ and $0.45$ and
with widths $0.2, 0.15, 0.08$ and $0.04$ respectively; and one type N
strip, where $\dot\theta/n\in[0.47, 0.5]$, $\dot\theta = 0.5 n$ being the
position of the first kink. The reasoning behind the choice of these values
is given in sections~\ref{reg_a} and~\ref{reg_n}. The expansion point
$\theta_e$ for the series solution is always the centre of the strip, and is shown by `{\sf x}'
in figure~\ref{bands}.

\subsection{Region H: HEM applies}
\label{reg_a}

The HEM, which is a fixed timestep implementation of the Frobenius method,
can in principle be used to solve equation~(\ref{ode}) in regions
of the phase plane, $(\theta, \dot\theta)$, where the functions $\atri$ and
$\atide$ are analytic. From its definition, $\atri(\theta, t)$ is an entire
function, so this is no bar to using the HEM. However, $\atide(\dot\theta)$
is $C^1$ in $\dot\theta$, with its second derivative being undefined at
$\dot\theta/n = 1/2, 1, 3/2,\ldots, 9/2$. Hence we approximate $\atide$ by its Taylor 
series of degree $D_{\rm TID}$ about a point $\dot\theta = \dot\theta_e$, bearing in 
mind that this will only be good for values of $\dot\theta$ far from the kinks. In the work 
reported here, our error criterion (see below) is satisfied with $D_{\rm TID} = 25$ for all strips.

With these provisos, the HEM works well. We briefly review the method here;
full details can be found in~\cite{hem}.

Let the state vector $\v{x}(t) = \left(\theta(t), \dot\theta(t)\right)$ and
define $t_i = ih,\; i\in\mathbb N_0$ where $h>0$ is a timestep. The ODE
allows us to compute, recursively, derivatives of $\theta(t)$ of all orders,
far from the kinks.
Let the $j$-th time derivative of $\theta(t)$ be written as $\theta^{(j)}(t)$.
We then write the degree-$D_s$ series solution to equation~(\ref{ode}) as
\begin{equation}
\v{x}(t_i) = \v{p}_i\left(\v{x}(t_{i-1})\right) =\v{x}(t_{i-1}) 
+ \sum_{j=1}^{D_s}\frac{h^j}{j!} \v{f}_j\left(\v{x}(t_{i-1}), t_{i-1}\right) + O\left(h^{D_s+1}\right)
\label{ppm}
\end{equation}
where $D_s > 1$ ($D_s = 1$ gives the Euler method), and
$\v{f}_j(\v{x}, t) = \left(\theta^{(j)}(t),\, \theta^{(j+1)}(t)\right)$.
Note that $\v{f}(\v{x}, t)$ depends on $t$ as well as $\v{x}$ because the ODE is
non-autonomous. For the same reason, $\v{p}$ depends on $i$. To satisfy 
the error criterion, the following values of $D_s$ are used in the various bands: 
\begin{table}[h!]
\centering
\begin{tabular}{lcccccccccc}
Band  & 0  & 1  & 2  & 3  & 4  & 5  & 6  & 7  & 8  & 9\\
$D_s$ & 16 & 15 & 15 & 14 & 14 & 14 & 15 & 16 & 17 & 17\\
\end{tabular}
\end{table}

Equation~(\ref{ppm}) allows us to advance the solution by a time $h$. Our 
immediate objective, however, is to compute the Poincar\'{e} map for the ODE in
an efficient way. Let us denote the period in $t$ of $\atri$ as 
$T_0 = 2\pi/n$ (which is one `Mercury year'), and put $\v{x}_k = \v{x}(k T_0)$.
Then the Poincar\'{e} map, $\v{P}$, is defined by $\v{x}_{k+1} = \v{P}(\v{x}_k)$.
Iteration of $\v{P}$ many times starting from a given initial condition
$\v{x}_0 = (\theta(0), \dot\theta(0))$ enables us to generate a sequence
of `snapshots' of the state variables as the system evolves, from which we
can deduce the attractor in which the system is eventually captured. For
details of capture detection, see below. Using a
large number $I$ of initial conditions, one can then estimate the
probability of capture in any of the possible steady-state solutions.

In order to satisfy the error criterion, timestep $h$ needs to be sufficiently small.
As with all the parameters mentioned so far, there is a trade-off between
speed and accuracy; a choice of $h = T_0/M$ with $M = 24$ is a good
compromise in practice. Finally, then, the map $\v{P}$ can be built up from the functions
$\v{p}_i$ by $\v{P}(\v{x}) = \v{p}_M\circ\v{p}_{M-1}\circ\ldots\circ\v{p}_1(\v{x})$,
and this is the HEM.

The expressions for $\v{p}_i$ are derived by computer algebra and are initially 
very large, but most of the terms are negligible and hence can be pruned away.
By `term' we mean here a polynomial in $\dot\theta$ --- these are then multiplied
by powers of $\sin 2\theta, \cos 2\theta$ in order to make up $\v{p}_i$.
In practice, every term whose magnitude is less than $10^{-16}$ at the largest and
smallest values of $\dot\theta$ within a strip is deleted. The value of
$10^{-16}$ is chosen because the usual double precision arithmetic is carried out to
around 16 significant figures. The resulting expression is
then converted to Horner form~\cite{Numrec} for efficient evaluation. Typically,
after pruning has been carried out, the expressions for $\v{p}_i$ contain 20 -- 35 terms
and are of the form: $\theta$-component $=\theta + r_3 + s r_5 + c r_5' + sc r_4 + s^2 r_2$,
$\dot\theta$-component $=r_4' + s r_6 + c r_6' + sc r_4'' + s^2 r_4''' +s^2c r_1 +s^3 r_1',$
where $r_k$, $r_k'$ etc.\ are polynomials in $\dot\theta(t_{i-1})$ of degree $k$, and
$s = \sin 2\theta$, $ c = \cos2\theta$.

We now describe the error criterion used. The final version of
the code for computing the Poincar\'{e} map via the HEM is compared to a
high-accuracy, standard numerical computation of the same thing. The full
expression for $\atide$ is used in the accurate numerical computation, not its series
approximation. The numerical algorithm used is the standard Runge-Kutta method
as implemented in the computer algebra
software Maple, computing to 25 significant figures and with absolute and
relative error tolerances of $10^{-15}$. The numerical and HEM computations
of $\v{P}(\v{x})$ are then compared, for each of the type H strips, using 250 uniformly
distributed random values of $\v{x}$ in each strip. The comparison gives
the maximum value of the modulus of the difference between each component,
computed both ways, over the 250 random points. The maximum difference observed
over all strips is assumed to be representative of the overall maximum
difference. Its values are about $3\times 10^{-14}$, $1.4\times 10^{-13}$ in the $\theta$- and
$\dot\theta$-components respectively.

\subsection{Region N: numerical method must be used}
\label{reg_n}

From figure~\ref{bands}, it can be seen that the HEM can be used for about 89\% of
$\dot\theta\in [0, 5n]$, but in the strips surrounding the kinks, $\theta \in [0, \pi)$,
$\dot\theta/n \in [(q+2)/2-0.03, (q+2)/2+0.03],\; q\in\Qtid$, the type N
strips, a purely numerical method has to be used.

It is possible that the ODE may be
stiff here, so we choose two numerical methods and compare the results. The
methods used are: (1) an explicit Runge-Kutta (RK) method due to Dormand and
Price, as described in~\cite{hnw}, and (2) the Adams method/Backward
Differentiation Formulae (BDF)~\cite{lsode}, with the ability to switch automatically between them.
The Adams method is an explicit predictor-corrector method, which, along
with RK, is suitable for non-stiff problems, whereas BDF is suitable for stiff problems.

In practice, with the parameters in Table~\ref{numvals}, BDF is rarely
needed, so the comparison between (1) and (2) above comes down to comparing
the RK and Adams methods. The implementation of Adams used~\cite{lsode} is approximately
1.6 times slower than RK for this problem, but the probabilities obtained from integration
starting from the same 3200 random points in $\mathcal Q$ are in good agreement --- see 
Table~\ref{rka} --- so, for type N strips, we choose RK, for which we fix
both the absolute and relative error tolerances to be $2\times 10^{-14}$. This
value is chosen to be comparable with the error entailed by polynomial
interpolation --- see below.

The probabilities obtained are not identical, neither should we expect them
to be. The long-term fate of a given trajectory depends very sensitively on
the details of its computation, implemented for very long integration
times. What is important in the end is the probabilities obtained, and
Table~\ref{rka} shows these to be robust against the algorithm used to
solve the ODE.

\begin{table}[h]
\centering
\begin{tabular}{|l|l|l||l|l|l|} \hline
\multicolumn{3}{|c||} {\raisebox{2.2ex}{ }Runge-Kutta} & \multicolumn{3}{c|} {\raisebox{2.2ex}{ }Adams} \\ \hline
\raisebox{2.7ex}{ }$\mtdn$ & Probability, \% & 95\% c.i. & 
$\mtdn$ & Probability, \% & 95\% c.i.\\ \hline
1/2 & 0.8125 & 0.3110 & 1/2 & 0.8750 & 0.3227\\
1   & 27.44  & 1.546  & 1   & 27.38 & 1.545\\
3/2 & 43.44  & 1.717  & 3/2 & 43.81 & 1.719\\
2   & 22.03  & 1.436  & 2   & 21.91 & 1.433\\
5/2 & 5.063  & 0.7596 & 5/2 & 5.094 & 0.7618\\
3   & 1.094  & 0.3604 & 3   & 0.750 & 0.2989\\
7/2 & 0.031  & (0.0612) & 7/2 & 0.094 & (0.1060)\\
4   & 0.094  & (0.1060) & 4   & 0.094 & (0.1060)\\ \hline
\end{tabular}
\caption{Comparison of probability of capture by the eight attractors with
$\mtdn = 1/2, 1, \ldots, 4$. These were computed
using the Runge-Kutta and the Adams method in N-type strips, with HEM being used elsewhere.
The same 3200 uniformly distributed random initial conditions were used in both cases.
Confidence intervals in brackets are unreliable since there are too few data points
in the case of these attractors.}
\label{rka}
\end{table}

Computation time can be saved by efficient calculation of $\atri$ and
$\atide$. The expression for $\atri$ can be converted into a polynomial of
degree 1 in $\cos 2\theta$ and $\sin 2\theta$, and degree 8 in $\cos nt/2$
and $\sin nt/2$. In Horner form, this polynomial can be evaluated efficiently,
using 2 sin/cos evaluations, 16 addition and 36 multiplication operations. 

Evaluation of $\atide(\dot\theta)$ in the obvious way is computationally expensive,
since it requires the calculation of one fractional power per term in
equation~(\ref{sec_ttid}) --- nine in all. A more efficient way to evaluate it is:
\begin{description}
\item[Case 1:] if $\dot\theta$ is far from a kink, then use a pre-computed
Chebyshev polynomial fit~\cite{Numrec} to $\atide$ --- the function is very
smooth here;
\item[Case 2:] if $\dot\theta$ is close to a kink, compute the contribution to $\atide$
from that kink exactly, according to the appropriate single term in the sum~(\ref{sec_ttid}),
with the effect of the remaining kinks being replaced by a Chebyshev polynomial fit.
\end{description}
Hence, at most one fractional power is computed per evaluation of $\atide$.
In practice, we use a polynomial only (Case 1 above), of degree 25, unless $2\dot\theta/n$ is 
within $\pm 0.08$ of an integer, when Case 2 applies. In Case 2, we use a polynomial of
degree 7 to fit the remaining terms. 
The resulting absolute error is no more than $4\times 10^{-14}$. For comparison, the ratio
of the CPU time taken to evaluate $\atide$ directly, and via polynomial
fitting, is about 5.1.

A comparison of the timings in different circumstances is given in Table~\ref{timings},
for which 1000 random initial conditions were used. The CPU time taken to
iterate $\v{P}$ $10^5$ times starting from each of these was measured, with
any data in which the trajectory moved from a type H strip to one of type
N, or vice versa, being rejected.

\begin{table}[htbp]
\centering
\begin{tabular}{|l|c|l|} \hline
\multicolumn{1}{|c|}{Method} & Strip & Mean time for $10^5$ iterations of $\v{P}$, CPU-sec\\ \hline
{\raisebox{2.5ex}{}Runge-Kutta, slow $\atide$} & N & \hskip 0.7in 61.1\\ \hline
{\raisebox{2.5ex}{}Runge-Kutta, fast $\atide$} & N & \hskip 0.7in 20.5\\ \hline
{\raisebox{2.5ex}{}HEM}			       & H & \hskip 0.7in 0.313\\ \hline
\end{tabular}
\caption{Timings for the computation of $10^5$ iterations of the Poincar\'{e} map,
$\v{P}$, in various circumstances. These are: type H or type N strips, and,
in type N strips,
with or without the fast computation of $\atide$ described in section~\ref{reg_n}.
The mean was taken over 895/105 random initial conditions for type H/N strips
respectively.}
\label{timings}
\end{table}

\subsection{Capture test}
\label{cap_crit}

\begin{figure}[!ht]
\centering 
\includegraphics*[width=4in]{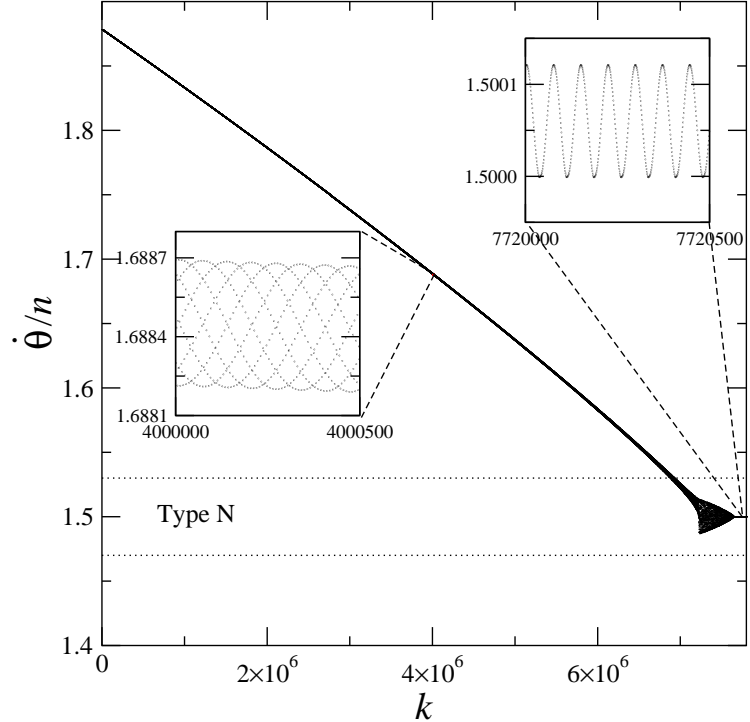}
\caption{An illustration of capture dynamics: $\dot\theta_k$, with initial
condition $\v{x}_0 = (0, 49) = (0, 1.878 n)$, is plotted against $k$, for $0\leq k\leq 7.8\times
10^6$. In this case, capture takes place in a solution
with $\mtdn \approx 3/2$. Two small regions are shown on a magnified scale.
On the left, we see the pre-capture dynamics and on the right,
post-capture. The HEM was used for all computations outside the strip $1.47
\leq\dot\theta/n\leq 1.53$, marked `Type N'.}
\label{capture}
\end{figure}

We now describe the test used to detect capture of a solution. 
Figure~\ref{capture} shows $\dot\theta_k = \dot\theta(kT_0)$ plotted against $k$
for $0\leq k\leq 7.8\times 10^6$. In this case, the initial condition was
$\v{x}_0 = (0, 49)$ and capture took place after about $7.7\times 10^6$
iterations of the Poincar\'{e} map, which corresponds to $7.7\times 10^6
T_0 = 1.85\times 10^6$ yr.

There are several ways that capture could be detected. We choose to divide the $\dot\theta_k$
dataset into blocks of length $L$ and compute the least squares gradient of each
block.  As can be seen from figure~\ref{capture}, this gradient will be
negative pre-capture, and close to zero post-capture. 

In detail the test is as follows. Let $\bar{y}_j$ be the mean value of
$\dot\theta_k/n$ over the $j$-th block of length $L$, and let $m_j$ be the
least squares gradient of $\dot\theta_k$ plotted against $k$ in block
$j$. Capture is deemed to have taken place when 
$$\left|2 \bar{y}_j - \left[2\bar{y}_j\right]\right| < \varepsilon_i
\;\;\;\mbox{ and }\;\;\;
\left|m_j\right| < \varepsilon_m
\;\;\;\mbox{for $K$ successive blocks}.$$
Here, $[x]$ is the nearest integer to $x$, and the factor of two occurs in
the first expression because capture can take place at integer or
half-integer values of $\dot\theta/n$. In practice, we choose $L = 10,000$,
$K = 8$, $\varepsilon_i = 10^{-3}$ and $\varepsilon_m = 3\times10^{-7}$.
The speed of a probability computation depends on the choice of these
parameters.

\begin{figure}[!ht]
\centering 
\includegraphics*[width=4in]{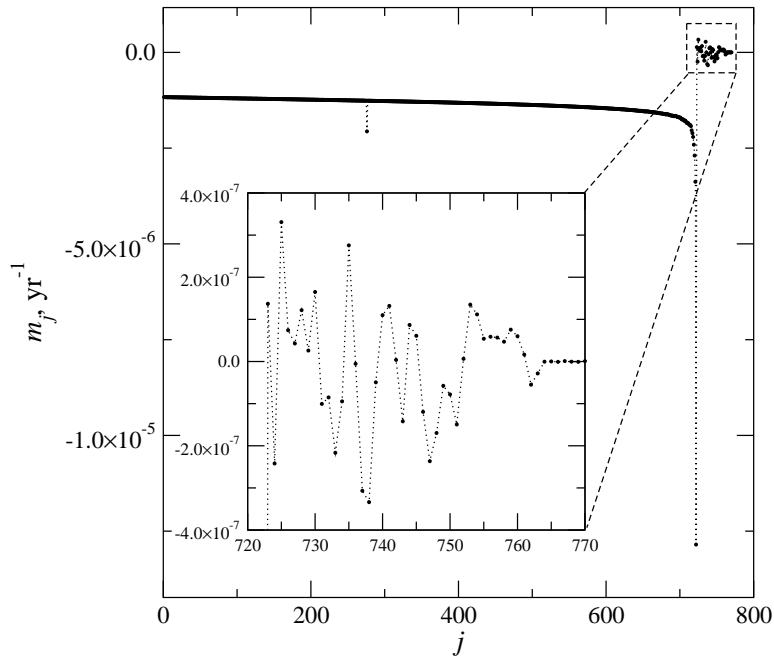}
\caption{An illustration of the capture test in practice.
The least-squares gradient of $\dot\theta_k$, $m_j$, versus the block number, $j$,
for the trajectory in figure~\ref{capture}. Blocks are of length 10,000
iterations of the Poincar\'{e} map and if one of the criteria for capture is
that $|m_j| < 3\times 10^{-8}$ for eight successive blocks, then this happens in
blocks 763--770 inclusive. Hence, capture in this case is detected after $7.7\times 10^6$
iterations of the Poincar\'{e} map.}
\label{lsm_capture}
\end{figure}

Figure~\ref{lsm_capture} shows a typical plot of $m_j$ versus block number $j$,
from $j = 0$ until capture, which takes place at $j = 770$.

\section{An illustrative probability of capture computation}

To illustrate how the foregoing works in practice, and to show some
timings, we now give results of a calculation of probability
of capture for the parameters in Table~\ref{numvals}.

As in~\cite{hem}, we use the CPU-sec as a unit of time, which is defined in
terms of the following sum:
$$S(m) = \sum_{i=1}^m\frac{(i+1)(i+3)}{i(i+2)(i+4)(i+6)},$$
whose evaluation requires $5m$ multiplications and $6m-1$ additions.
We define 1 CPU-sec to be the CPU time taken to evaluate $S(N_c)$, where $N_c = 5.96\times 10^7$,
which is the CPU time taken to do this computation on the computer used to do
some of the calculations in this paper. The time taken can vary according
to circumstances, e.g. the loading, the type of processor
and so on. Hence, care has to be taken in the codes to
scale the CPU-sec appropriately for the particular hardware used: the
computation of $S(N_c)$ is timed after every successful capture, and the
CPU-sec scaling factor is updated on each occasion.

One other point to note is that a computation of probability is trivial to
parallelise. The initial conditions in $\mathcal Q$ are generated by a
pseudo-random number generator, and a different sequence of pseudo-random
numbers can be produced just by changing the seed. Hence, by running $N$ copies of the
program on $N$ separate processors with $N$ different seeds, $N$ times the
number of initial conditions can be investigated at the same time.

The CPU time taken to integrate until capture depends strongly on the initial
condition, and in particular, the proportion of the integration that
is carried out in type N strips (which is a slow process) as opposed to type H
strips (where it is fast). By considering 57,600 random initial conditions, we
find the following:
$$\mbox{Mean time to capture: 1156 CPU-sec, with standard deviation 1190 CPU-sec}$$
$$\mbox{Maximum time to capture: 5161 CPU-sec; minimum time to capture: 1.188 CPU-sec}.$$
$$\mbox{Total iterations in type N strips: $\sim 1.6\times 10^{11}$;
total iterations in type H strips: $\sim 1.2\times 10^{12}$}.$$
From the above, we see that overall, about $1.4\times 10^{12}$ iterations of the
Poincar\'{e} map were needed to estimate the probability of capture for
57,600 initial conditions, and that 12\% of these were in type~N strips,
with the remaining 88\% being in type~H strips. These values have a
sensitive dependence on the capture parameters.

We can now estimate the factor by which our approach speeds up a typical
probability of capture computation. Let $T_{\mbox{mix}}(n)$ be the CPU time
taken to iterate the Poincar\'{e} map $n$ times using HEM in type H strips
and the chosen numerical method in type N strips, and let $T_{\mbox{num}}(n)$ 
be the time taken when using the numerical method everywhere. Then, using the data
above and from Table~\ref{timings}, we estimate that 
$$\frac{T_{\mbox{mix}}(n)}{T_{\mbox{num}}(n)} \approx \frac{n (0.12\times 20.5 +
0.88\times 0.313)10^{-5}}{20.5 n\times 10^{-5}}\approx \frac{2}{15}.$$
Our approach therefore speeds up this computation by a factor of about $7.5$.

\begin{figure}[!ht]
\centering 
\includegraphics*[width=3in]{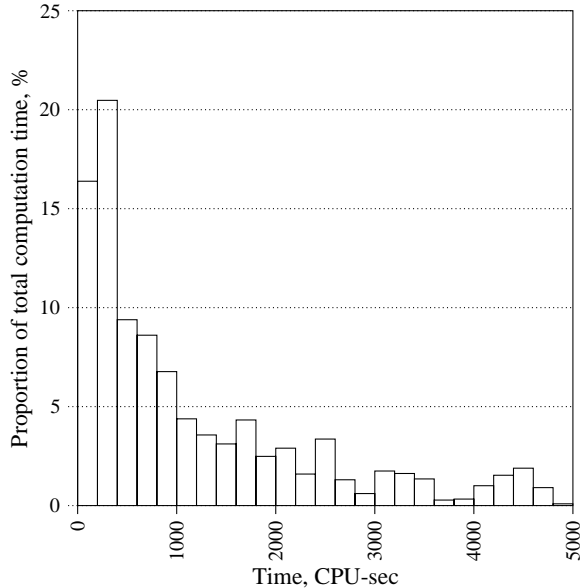}
\caption{Histogram showing the probability that the computation time for iteration to capture
lies in various ranges. The interpretation is that, for instance, a CPU time
of between 0 and 200 CPU-sec is required to establish the capture of about 16.3\%
of trajectories (leftmost bar). Four values of CPU-time $> 5000$ CPU-sec were 
excluded from the histogram and so the dataset contains 57~596 points.
The modal value of CPU-time is in the range 200--400 CPU-sec.}
\label{time_hist}
\end{figure}

A histogram of capture times is given in figure~\ref{time_hist}. Let $\tau\left(\v{x}_0\right)$
be the number of CPU-sec required to iterate the Poincar\'{e} map, starting
from $\v{x}_0$, until the capture criterion of section~\ref{cap_crit} is met. In
order to produce Figure~\ref{time_hist}, 57,600 random values of $\v{x}_0\in\mathcal Q$
were generated and $\tau\left(\v{x}_0\right)$ was computed for each. This figure
is a histogram of the proportions of the values of $\tau\left(\v{x}_0\right)$ that lie
in the ranges 0--200, 200--400, \ldots, 4800--5000 CPU-sec.

Finally, we give an estimate of the probabilities themselves in Table~\ref{prodrun}. The
95\%-confidence interval~\cite{Wal} is defined such that the probability of capture in a
given attractor $A$ lies in the interval $[\hat{p} - \Delta p, \hat{p} + \Delta p]$, with
95\% confidence. Here,
$\Delta p = 1.96\sqrt{\hat{p}(1-\hat{p})/I}$, in which $I = 5.76\times 10^4$ is the number of
initial conditions and $\hat{p}$ is the proportion of these initial conditions that end up in $A$.
For these values, in the worst case, which is $\hat{p} = 1/2$, $\Delta p \approx 0.41$\%.

\begin{table}[h]
\centering
\begin{tabular}{|l||l|l|l|} \hline
\multicolumn{1}{|c||} {\raisebox{2.2ex}{ }$p/q$} & \multicolumn{1}{c|} {\raisebox{2.2ex}{ }Number}
& Probability, $\hat p$, \% & 95\% c.i. \\ \hline
1/2	& 349	& 0.6059	& 0.063\\
2	& 16292	& 28.285	& 0.368\\
3/2	& 24344	& 42.264	& 0.403\\
4	& 13252	& 23.007	& 0.344\\
5/2	& 2775	& 4.8177	& 0.175\\
6	& 473	& 0.8212	& 0.074\\
7/2	& 89	& 0.1545	& 0.032\\
8	& 26	& 0.0451	& 0.017\\
9/2	& 0	& 0.0000	& 0.000\\ \hline
\end{tabular}
\caption{Probability of capture by various attractors with
$\mtdn = 1/2, 1, \ldots, 4$; as expected, 9/2 was never observed. These were computed
using the Runge-Kutta method~\cite{hnw} in N-type strips, with HEM being used elsewhere.
We used 57,600 uniformly distributed random initial conditions in
$\mathcal Q = [0, \pi]\times[0, 5n]$.}
\label{prodrun}
\end{table}

\section{Conclusions}

It was not the main purpose of this paper to give results of the computation of
probabilities of capture for various sets of parameters --- rather, we wanted
to show how this computation, in the case of constant eccentricity $e$, can
be accelerated. This requires us to solve the ODE~(\ref{ode}),
and bearing in mind that this is periodic in time with period $T_0 = 2\pi/n$, our
objective is to compute the Poincar\'{e} map, which advances the state
variables $\v{x} = (\theta, \dot\theta)$ by a time $T_0$, as fast as possible.
This enables us to generate a sequence of values $\v{x}_k$, from which salient
information about the dynamics can be deduced.

The dissipation term, $\atide(\dot\theta)$, which is
a $C^1$ function of the angular velocity $\dot\theta$, and additionally varies 
rapidly in the ranges of $\dot\theta$ around the so-called `kinks', complicates
the computation, and requires the use of a standard numerical ODE solver for some
of the time.  We point out ways in which these purely numerical computations can
be carried out efficiently, by streamlining the calculation of the triaxiality and tidal
accelerations. However, in about 88\% of the phase plane the high-order Euler 
method~\cite{hem} can be used, and this speeds up the computation significantly.

The present case should be contrasted with the constant time lag model,
in which dissipation is just proportional to $\dot\theta - \omega$ with $\omega$ a
constant. In the light of its simplicity, this has been used in many publications, 
for instance, \cite{hem} and references therein. In that case, only a single set of
integrators $\v{p}_i$, defined in equation~(\ref{ppm}), was needed to build up the 
Poincar\'{e} map. In other words, there was only one strip, which was the whole of 
$\mathcal Q$, and the HEM could be used everywhere.

Compare that with the current case using the parameters in Table~\ref{numvals}.
Here, $\mathcal Q$ has to be divided into 48 strips and an integrator set up for each.
Setting up the codes for the HEM in each strip, which is done by computer
algebra, itself takes significant time, but the pay-off is an
increase in speed by a factor of approximately 65 in the strips where HEM
can be used, compared to using a standard numerical method.

A probability of capture computation requires a capture detection
algorithm, and one has been described. It is based on the fact that after
capture, the angular velocity has no underlying decay: a captured solution
has a close to constant mean value of $\dot\theta_k = \dot\theta(k T_0)$ for all
$k$ greater than the value at which capture takes place.
Typically, the mean is taken over $10^4$ successive $k$ values.

A test run of our algorithm, in which 57,600 initial conditions were
iterated until capture took place, reveals that the overall speed of computation
is faster by a factor of about $7.5$ compared to using a standard numerical
algorithm alone.

An important question left unanswered in this paper is `What is the nature
of the solutions in which capture takes place?' The answer turns out to be
more complicated than expected: periodic solutions with mean
$\dot\theta/n \approx -1, -1/2, 1/2,\ldots, 4$ have been computed (high-accuracy
numerics are needed). There is also numerical evidence for the existence of attracting
solutions whose period is not a
small integer multiple of $T_0$. The dynamics of solutions to equation~(\ref{ode}),
both pre- and post-capture, will be subject of a future publication~\cite{paper2}.

\end{document}